\documentclass[a4paper,12pt]{amsart}
\usepackage{mathtools}

\title
{Convex caterpillars are Schur-positive}

\author{%
	Yuval H. Khachatryan-Raziel%
}

\thanks{Supported by Israel Science Foundation grant no. 1970/18}

\usepackage[utf8]{inputenc}
\usepackage[english]{babel}
\usepackage{amssymb,amsmath,amsthm,amsfonts, tikz}

\makeatletter
\DeclareRobustCommand*{\bfseries}{%
	\not@math@alphabet\bfseries\mathbf
	\fontseries\bfdefault\selectfont
	\boldmath
}
\makeatother

\theoremstyle{plain}

\newtheorem{theorem}{Theorem}[section]

\newtheorem{prop}[theorem]{Proposition}
\newtheorem{lemma}[theorem]{Lemma}
\newtheorem{corollary}[theorem]{Corollary}

\newtheorem{problem}[theorem]{Problem}
\theoremstyle{definition}
\newtheorem{definition}[theorem]{Definition}
\newtheorem{defn}[theorem]{Definition}
\newtheorem{example}[theorem]{Example}

\newtheorem{observation}[theorem]{Observation}

\newtheorem{fact}[theorem]{Fact}

\numberwithin{equation}{section}

\newcommand{\Q}{\mathcal Q}
\newcommand{\F}{\mathcal F}
\newcommand\x{{\mathbf x}}
\newcommand{\bbq}{\mathbb{Q}}

\DeclareMathOperator\Des{Des}

\DeclareMathOperator\SYT{SYT}

\newcommand{\NN}{\mathbb{N}}

\newcommand{\BBB}{{\mathcal{B}}}

\newcommand{\ch}{\operatorname{ch}}

\newcommand{\symm}{{\mathfrak{S}}}
\newcommand{\Ct}{{Ct}}

\begin{document}
	
\maketitle

\begin{abstract}
	A remarkable result of Stanley shows that the set  of maximal chains in the non-crossing partition lattice of type $A$ is Schur-positive, where descents are defined by a distinguished edge labeling. A bijection between these chains and labeled trees was presented by Goulden and Yong.
	Using Adin-Roichman's variant of Bj\"orner's $EL$-labeling, we show that the subset of maximal chains in the non-crossing partition lattice of type $A$, whose underlying tree is a convex caterpillar, 
	is Schur-positive.
\end{abstract}


\tableofcontents
	
	\section{Introduction}

	A symmetric function is called {\em Schur-positive} if all the
	coefficients in its expansion in the basis of Schur functions are nonnegative.
	Determining whether a given symmetric function is
	Schur-positive is a major problem in contemporary algebraic
	combinatorics~\cite{Stanley_problems}. 	
	
	With a set $A$ of combinatorial objects, equipped with a {\em descent map} $\Des: A \to 2^{[n-1]}$, one associates the
	quasi-symmetric function
	\[
	\Q(A) := \sum\limits_{\pi\in A} \F_{n,\Des(\pi)}
	\]
	where 
	$\F_{n,D}$ (for $D \subseteq [n-1]$) are Gessel's {\em fundamental quasi-symmetric functions;} see Subsection~\ref{sec:prel_quasi} for more details. 
	The following problem  is long-standing. 
	
	\begin{problem}\label{prob:symmetric}
		Given a set $A$, equipped with a descent map, 
		is $\Q(A)$ symmetric?
		In case of an affirmative answer, is it Schur-positive?
	\end{problem}

	
	
	\medskip
	
	Of special interest are Schur-positive sets of maximal chains. Maximal chains in a labeled poset $P$ are equipped with a natural descent map.
	A well-known conjecture of Stanley~\cite[III, Ch. 21]{Stanley_thesis}
	implies that all examples of Schur-positive labeled posets in this sense
	correspond to intervals in the Young lattice. 

	
	Another way to equip the set of maximal chains with a descent map is using a labeling of the edges in the Hasse diagram.
	A classical example of a Schur-positive set 
	of this type, 
	the set of all maximal chains in the non-crossing partition lattice of type $A$,
	was given by Stanley~\cite{StanleyPFnNCPL}. 
	An $EL$ edge-labeling of this poset was presented in an earlier work of Bj\"orner~\cite{Bjorner}; see also~\cite{BW, Mcnamara, AR_OMCNCPL}.


	
	
	

	\medskip
	
	The goal of this paper is to present an interesting set of maximal chains in the non-crossing partition lattice $NC_n$ 
	(equivalently: a set of edge-labeled trees) 
	which is Schur-positive. 
	We will use a variant of Bj\"orner's $EL$-labeling, presented in~\cite{AR_OMCNCPL}.
	
	\medskip
	
	It is well known that 
	maximal chains in the non-crossing partition lattice may be interpreted as factorizations of the $n$-cycle $(1,2,\dots,n)$
	into a product of $n-1$ transpositions.
	
	\begin{definition}
		A factorization $t_1\cdots t_{n-1}$ of the $n$-cycle $(1,2,\dots,n)$ as a product of transpositions is called {\em linearly ordered} if, 
		for every $1 \le i \le n-2$, $t_i$ and $t_{i+1}$ have a common letter.
	\end{definition}
	
	This definition is motivated by Theorem~\ref{thm:linear} below.
	Denote the set of linearly ordered factorizations of  $(1,2,\dots,n)$ by $U_n$.

	\begin{prop}
		For every $n\ge 1$, the number of linearly ordered factorizations of the $n$-cycle $(1,2,\dots,n)$ is 
		\[
		|U_n|=n2^{n-3}.
		\]	\end{prop}
	Our main result is
	
	\begin{theorem}\label{thm:main}
		The set of linearly ordered factorizations of the $n$-cycle $(1,2,\dots,n)$
		satisfies 
		$$\Q(U_n)=\sum_{k=0}^{n-1} (k+1)s_{(n-k,1^k)},$$
		where the descent set of any $u\in U_n$ is defined by the edge labeling of~\cite{AR_OMCNCPL}.
		In particular, $U_n$ is Schur-Positive.
	\end{theorem}

	
	It should be noted that Theorem~\ref{thm:main} does not follow from Stanley's proof of the Schur-positivity of the set of all maximal chains in $NC_n$. In fact, Stanley's action on maximal chains does not preserve linearly ordered chains. 
	
	We prove Theorem~\ref{thm:main}, by 
	translating it into the language of geometric trees called convex caterpillars.
	


	\begin{definition}
		A tree is called a {\em caterpillar} if the subgraph obtained by removing all its leaves is a path. This path is called the \emph{spine} of the caterpillar.
	\end{definition}
	
	\begin{definition}
		A {\em convex caterpillar} of order $n$ is a caterpillar drawn in the plane such that
		\begin{enumerate}
		\item[(a)]
		the vertices are in convex position (say, the vertices of a regular polygon) and labeled ${1,\ldots,n}$ clockwise;
		\item[(b)]
		the edges are drawn as non-crossing straight line segments; and
		\item[(c)]
		the spine forms a cyclic interval $(a,a+1),(a+1,a+2),\ldots,(b-1,b)$ in $[n]$.
		\end{enumerate}
		
		Denote by $\Ct_n$ the set of convex caterpillars of order $n$.
	\end{definition}
	
\begin{figure}[hbt]
		\begin{center}
		\begin{tikzpicture}[scale=1]
		\fill (2.4,3.4) circle (0.1) node[above]{\tiny 1};
		\fill (3.4,2.4)	circle (0.1) node[right]{\tiny 2};
		\fill (3.4,1) circle (0.1) node[right]{\tiny 3}; 
		\fill (2.4,0) circle (0.1) node[below]{\tiny 4};
		\fill (1,0) circle (0.1) node[below]{\tiny 5};
		\fill (0,1)	circle (0.1) node[left]{\tiny 6};
		\fill (0,2.4) circle (0.1) node[left]{\tiny 7};
		\fill (1,3.4) circle (0.1) node[above]{\tiny 8};
			
		\draw (3.4,2.4)--(3.4,1);
		\draw (0,1)--(1,3.4)--(1,0);
		\draw (3.4,2.4)--(2.4,0); 
		\draw (0,2.4)--(1,3.4);
	
		\draw[blue] (1,3.4)--(2.4,3.4); 
		\draw[blue] (3.4,2.4)--(2.4,3.4);
		\end{tikzpicture}
		\end{center}
		\caption{A convex caterpillar and its spine}\label{fig:1}
	\end{figure}
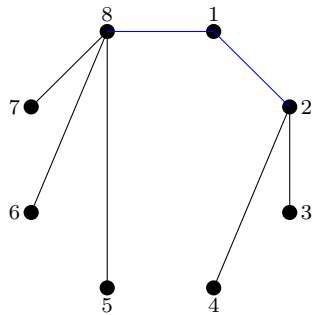

	\begin{example}
	Figure~\ref{fig:1} shows a convex caterpillar $c \in \Ct_8$, with spine consisting of the edges $(8,1)$ and $(1,2)$, forming a cyclic interval.
	\end{example}
	
	
	Goulden and Yong~\cite{Goulden-Yong} introduced a mapping from factorizations of $(1,2,\ldots,n)$ to non-crossing geometric trees.
	This mapping is not injective: in order to recover the factorization from the tree, one has to choose a linear extension of a certain partial order on the edges, which we call the {\em Goulden-Yong partial order}; see Definition~\ref{defn:GLOrder} below.
	
	In a previous work~\cite{YK-CenterPaper} we proved that the Goulden-Yong order is linear if and only if the geometric tree is a convex caterpillar; see Theorem~\ref{thm:linear} below.
	It follows that the Goulden-Yong map, restricted to the set $U_n$ of linearly ordered factorizations, is a bijection onto the set $\Ct_n$ of convex caterpillars of order $n$.
	
	\begin{definition}\label{defn:Des-c}
		The {\em descent set} of a linearly ordered factorization $u = (t_1, \ldots, t_{n-1}) \in U_n$ is
		\[
		\Des(u) := \{i\in [n-2] \,:\, t_i=(b,c) \text{ and } t_{i+1}=(b,a) \text{ with }\ c > a\}.
		\]
	\end{definition}

	\begin{example} 
	    The convex caterpillar $c \in \Ct_8$, drawn in Figure~\ref{fig:1},
		corresponds to the linearly ordered word
		\[
		u = \left( (7,8),(6,8),(5,8),(1,8),(1,2),(2,4),(2,3) \right) \in U_8,
		\]
		for which $\Des(u) = \{1,2,3,4,6\}$.
	\end{example}
	


	In \cite{AR_OMCNCPL}, the authors define a map $\phi$ from the set denoted here $U_n$ to the symmetric group $\symm_{n-1}$; for a detailed description see Subsection~\ref{sec:phi} below. 
	The map $\phi$ is an $EL$-labeling of the non-crossing partition lattice. This property, relations to Bj\"orner's $EL$-labeling 
	and other positivity phenomena  will be discussed in another paper.
	
	It turns out that our Definition~\ref{defn:Des-c} above fits nicely with this map.
	
	
	
	\begin{lemma}
		For any $u \in U_n$, 
		\[
	    \Des(\phi(u)) = \Des(u).
		\]
	\end{lemma}
	
	See Proposition~\ref{prop:phi_caterpillar_des} below.
	We further show that the number of caterpillars with a given descent set depends only on the cardinality of the descent set.	
	
	\begin{lemma}
		For every subset $J\subseteq [n-2]$,
		\[
		|\{c \in \Ct_n:\ \Des(c)=J\}|= |J| + 1.
		\]
	\end{lemma}

	These two key lemmas are  used to prove Theorem~\ref{thm:main}.
	

	\section{Background}
	In this section we provide the necessary definitions and historical background to explain the main results. More information can be found in the references.
	
	
	\subsection{Compositions, partitions and tableaux}
	\begin{definition}
		A {\em weak composition} of $n$ is a sequence $\alpha = \left(\alpha_1,\alpha_2,\dots\right)$ of non-negative integers such that $\sum_{k=1}^\infty\alpha_k = n$.
	\end{definition}
	\begin{definition}
		A {\em partition} of $n$ is a weakly decreasing sequence of non-negative integers $\lambda = (\lambda_1,\lambda_2,\dots)$ such that $\sum_{k=1}^\infty\lambda_k = n$. We denote $\lambda \vdash n$.
	\end{definition}
	\begin{definition}
    		The {\em length} of a partition $\lambda = \left(\lambda_1,\lambda_2,\dots \right)$ is the number of non-zero parts $\lambda_i$.
	\end{definition}

	

	
	For a skew shape $\lambda/\mu$, let $\SYT(\lambda/\mu)$ be the set
	of standard Young tableaux of shape $\lambda/\mu$. We use the
	English convention, according to which row indices increase from top to bottom
	(see, e.g., \cite[Ch.\ 2.5]{Sagan}). The {\em height} of a standard Young tableau $T$ is the number of rows in $T$.
	%
	The {\em descent set} of $T$ 
	is
	\[
	\Des(T):=\{i \,:\, i+1 \text{ appears in a lower row of } T \text{ than } i\}.
	\]

	

	
	
	\subsection{Symmetric and quasi-symmetric functions}\label{sec:prel_quasi}
	
	
	
	
	Let ${\bf x} := (x_1,x_2,\ldots)$ be an infinite sequence of commuting indeterminates. 
	Symmetric and quasi-symmetric functions in ${\bf x}$ can be defined over various (commutative) rings of coefficients, including the ring of integers; for simplicity we define it over the field $\bbq$ of rational numbers. 
	
	\begin{definition}                                 
		A {\em symmetric function} in the variables $x_1,x_2,\ldots$ is a formal power series $f({\bf x})\in \bbq[[{\bf x}]]$, of bounded degree, 
		such that for any three sequences (of the same length $k$) of positive integers, $(a_1,\ldots,a_k)$, $(i_1,\dots,i_k)$ and $(j_1,\dots, j_k)$, the coefficients of $x_{i_1}^{a_1}\cdots x_{i_k}^{a_k}$ and of $x_{j_1}^{a_1}\cdots x_{j_k}^{a_k}$ in $f$ are the same:
		\[
		[x_{i_1}^{a_1}\cdots x_{i_k}^{a_k}]f = 
		[x_{j_1}^{a_1}\cdots x_{j_k}^{a_k}]f.
		\]
	\end{definition}
	
	Schur functions, indexed by partitions of $n$, form a distinguished basis for $\Lambda^n$, the vector space of symmetric functions which are homogeneous of degree $n$;
	see, e.g., \cite[Corollary 7.10.6]{EC2}.
	A symmetric function in $\Lambda^n$ is {\em Schur-positive} if all the coefficients in its expansion in the basis $\{s_{\lambda} \,:\, \lambda \vdash n\}$ of Schur functions are non-negative.
	
	The following definition of a quasi-symmetric function can be found in~\cite[7.19]{EC2}.
	\begin{definition}                                 
		A {\em quasi-symmetric function} in the variables $x_1,x_2,\ldots$  is a formal power series $f({\bf x})\in \bbq[[{\bf x}]]$, of bounded degree, 
		such that for any three sequences (of the same length $k$) of positive integers, $(a_1,\ldots,a_k)$, $(i_1,\dots,i_k)$ and $(j_1,\dots, j_k)$, where the last two are {\em increasing}, the coefficients of $x_{i_1}^{a_1}\cdots x_{i_k}^{a_k}$ and of $x_{j_1}^{a_1}\cdots x_{j_k}^{a_k}$ in $f$ are the same:
		\[
		[x_{i_1}^{a_1}\cdots x_{i_k}^{a_k}]f = 
		[x_{j_1}^{a_1}\cdots x_{j_k}^{a_k}]f
		\]
		whenever $i_1 < \ldots < i_k$ and $j_1 < \ldots < j_k$. 
	\end{definition}
	Clearly, every symmetric function is quasi-symmetric, but not conversely: $\sum_{i<j}{x_i^2 x_j}$, for example, is quasi-symmetric but not symmetric.
	
	For each subset $D \subseteq [n-1]$ define the {\em
		fundamental quasi-symmetric function}
	\[
	\F_{n,D}(\x) := \sum_{\substack{i_1\le i_2 \le \ldots \le i_n \\ {i_j < i_{j+1} \text{ if } j \in D}}} 
	x_{i_1} x_{i_2} \cdots x_{i_n}.
	\]
	
	Let $\BBB$ be a set of combinatorial objects, equipped with a {\em descent map} $\Des: \BBB \to 2^{[n-1]}$ which associates to each
	element $b\in \BBB$ a subset $\Des(b) \subseteq [n-1]$. Define the
	quasi-symmetric function
	\[
	\Q(\BBB) := \sum\limits_{b\in \BBB} 
	\F_{n,\Des(b)}.
	\]
	With some abuse of terminology, we say that $\BBB$ is Schur-positive when $\Q(\BBB)$ is.
	
	\medskip
	
	The following key theorem is due to Gessel.
	
	\begin{theorem}{\rm \cite[Theorem 7.19.7]{EC2}}\label{G1} 
		For every shape $\lambda \vdash n$,
		\[
		\Q({\SYT(\lambda)})=s_{\lambda}.
		\]
	\end{theorem}
	
	
	\begin{corollary}\label{thm:Sp-combin}
		A set $\BBB$, equipped with a descent map $\Des: \BBB \to 2^{[n-1]}$, is Schur-positive if and only if there exist nonnegative integers $(m_{\lambda,\BBB})_{\lambda\vdash n}$ such that
		\begin{equation}\label{eq:equidist}
		\sum_{b \in \BBB} {\bf x}^{\Des(b)} 
		=\sum_{\lambda \vdash n} m_{\lambda,\BBB} \sum_{T\in \SYT(\lambda)} {\bf x}^{\Des(T)}. 
		\end{equation}
	\end{corollary}
	
	There is a dictionary relating symmetric functions to characters of the symmetric group $\symm_n$. The irreducible characters of $\symm_n$ are indexed by partitions $\lambda \vdash n$ and denoted $\chi^\lambda$. 
	The {\em Frobenius characteristic map} $\ch$ from class functions on $\symm_n$ to symmetric functions is defined by $\ch(\chi^{\lambda}) = s_{\lambda}$, and extended by linearity.
	Theorem~\ref{G1} may then be restated as follows:
	\[
	\ch(\chi^\lambda) = \sum_{T \in SYT(\lambda)} \F_{n,\Des(T)}.
	\]
	
	

	\subsection{Maximal chains in the non-crossing partition lattice}
	
	The systematic study of noncrossing partitions began with Kreweras \cite{Kreweras-NonCrossing} and Poupard \cite{Poupard-NonCrossing}.
	Surveys of results and connections with various areas of mathematics can be found in 
	\cite{Simion-NonCrossing} and \cite{Armstrong}.

    A {\em noncrossing partition} of the set $[n]$ is a partition $\pi$ of $[n]$ into nonempty blocks with the following property: 
	for every $a<b<c<d$ in $[n]$, if some block $B$ of $\pi$ contains $a$ and $c$ and some block $B^\prime$ of $\pi$ contains $b$ and $d$, then $B=B^\prime$. 
	Let $NC_n$ be the set of all noncrossing partitions of $[n]$.
	Define a partial order on $NC_n$, by refinement: 
	$\pi \leq \sigma$ if every block of $\pi$ is contained in a block of $\sigma$.
	This turns $NC_n$ into a graded lattice.
	
	An \emph{edge labeling} of a poset $P$ is function from the edges of the Hasse diagram of $P$ to the set of integers. Several different edge labelings of $NC_n$ were defined and studied 
	by Bj\"orner \cite{Bjorner}, Stanley \cite{StanleyPFnNCPL}, and Adin and Roichman \cite{AR_OMCNCPL}.
    Let $\Lambda$ be an edge labeling of $NC_{n+1}$,
    and let $F_{n+1}$ be the set of maximal chains in $NC_{n+1}$.
    For each maximal chain $\mathfrak{m} : \pi_0 < \pi_1 < \dots < \pi_n$ define
	\[
	\Lambda^*(\mathfrak{m}) := \left( \Lambda(\pi_0, \pi_1), \ldots, \Lambda(\pi_{n-1}, \pi_n) \right) \in \NN^n,
	\]
	with a corresponding {\em descent set}
	\[
	\Des(\mathfrak{m}) := \left\{ i \in [n-1] \,:\, \Lambda(\pi_{i-1}, \pi_i) >
	\Lambda(\pi_i, \pi_{i+1}) \right\}.
	\]

\medskip



	



	
	The noncrossing partition lattice is is intimately related to cycle factorizations.
	The $n$-cycle $(1,2,\ldots,n)$ can be written as a product of $n-1$ transpositions. 	There is a well known bijection between such factorizations and the maximal chains in $NC_{n+1}$; see, for example, \cite[Lemma 4.3]{NC_SurprisingLocations}. 
		A classical result of Hurwitz states that the number of such factorizations is $n^{n-2}$ \cite{HurwitzClassical, HurwitzReproof},
		thus equal to the number of labeled trees of order $n$.
		In the next section we will describe a connection between maximal chains and 	
	geometric trees.
	
	

		\section{The Goulden-Yong partial order}
	
\begin{figure}[hbt]
	\begin{center}
		\begin{tikzpicture}[scale=1]
		\fill (2.4,3.4) circle (0.1) node[above]{\tiny 1};
		\fill (3.4,1.7)	circle (0.1) node[right]{\tiny 2};
		\fill (2.4,0) circle (0.1) node[below]{\tiny 3};
		\fill (1,0) circle (0.1) node[below]{\tiny 4};
		\fill (0,1.7)	circle (0.1) node[left]{\tiny 5};
		\fill (1,3.4) circle (0.1) node[above]{\tiny 6};
			
		\draw (2.4, 3.4) -- (3.4, 1.7);
		\draw (3.4, 1.7) -- (2.4, 0);
		\draw (2.4, 3.4) -- (1, 0);
		\draw (1, 0) -- (0, 1.7);
		\draw (1, 0) -- (1, 3.4);
		\end{tikzpicture}
	\end{center}
	\caption{$G(w)$ for $w = ((1,4),(4,6),(4,5),(1,2),(2,3))$}\label{fig:2}
	\end{figure}
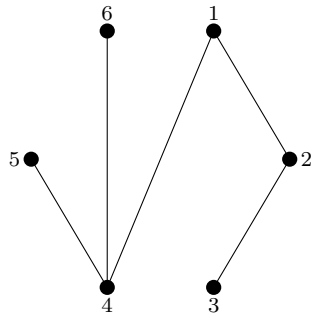

	   	With each sequence of $n-1$ different transpositions $w = (t_1,\dots,t_{n-1})$, associate a geometric graph $G(w)$ as follows.
   	The vertex set is the set of vertices of a regular $n$-gon, labeled clockwise $1,2,\dots,n$.
   	The edges correspond to the given transpositions $t_1,\dots,t_{n-1}$, where the edge corresponding to a transposition $t_k = (i,j)$ is the line segment connecting vertices $i$ and $j$.
   	See Figure~\ref{fig:2} for the geometric graph $G(w)$ corresponding to $w = ((1,4),(4,6),(4,5),(1,2),(2,3))$.
	
	Let $F_n$ be the set of all factorizations of the $n$-cycle $(1,2,\ldots,n)$ into a product of $n-1$ transpositions.
	Write each element of $F_n$ as a sequence $(t_1, \ldots,t_{n-1})$, where $t_1 \cdots t_{n-1} = (1,2,\ldots,n)$.
	The following theorem of Goulden and Yong gives necessary and sufficient conditions for a sequence of $n-1$ transpositions to belongs to $F_n$.
	
	\begin{theorem}\cite[Theorem 2.2]{Goulden-Yong}\label{thm:GY}
		A sequence of transpositions $w = (t_1, \ldots, t_{n-1})$ belongs to $F_n$ if and only if the following three conditions hold:
		\begin{enumerate}
			\item $G(w)$ is a tree.
			\item $G(w)$ is non-crossing,  namely: two edges may intersect only in common vertex.
			\item Cyclically decreasing neighbors: For every 
			$1 \leq i < j \leq n-1$, if $t_i = (a, c)$ and $t_j = (a, b)$ then $c >_a b$. Here $<_a$ is the linear order $a <_a a + 1 <_a \dots <_a a-1$.
		\end{enumerate}
	\end{theorem}

	For example, the graph in Figure~\ref{fig:2} corresponds to a sequence $w \in F_6$, and indeed satisfies the conditions of Theorem~\ref{thm:GY}.
	
	Note that a sequence $w = (t_1, \ldots, t_{n-1}) \in F_n$ carries more information than its Goulden-Yong tree $G(w)$: it actually defines a {\em linear order} on the edges, with the edge corresponding to $t_i$ preceding the edge corresponding to $t_j$ whenever $i < j$. 
	How much of that information can be retrieved from the tree?

	
	\begin{defn}\label{defn:GLOrder}
	Let $T$ be a non-crossing geometric tree 
    (namely, satisfying conditions 1 and 2 of Theorem~\ref{thm:GY}) 
	on the set of vertices of a regular $n$-gon, labeled clockwise $1,2,\dots,n$.
	Define a relation $\le_T$ on the set of edges of $T$ as follows:
	$(a,b) \le_T (c,d)$ if there exists a sequence of edges $(a,b) = t_0, \ldots, t_k =(c,d)$ $(k \ge 0)$ such that for every $0 \leq i \leq k-1$, $t_i = (x, z)$ and $t_{i+1} = (x,y)$ have a common vertex $x$ and $z >_x y$ as in condition 3 of Theorem~\ref{thm:GY}.
	\end{defn}
	
	\begin{lemma}\label{lemma:Goulden_Yong_Order_On_T}
	$\le_T$ is a partial order on the set of edges of $T$.
	\end{lemma}
	
	We use the following well-known fact to prove the statement.
	\begin{fact}\label{fact:TransitiveClosureProperty}
    Let $R$ be an anti-symmetric relation on a set $S$ such that for every $x,y\in S$ there is at most 
    one finite sequence $x = a_0,\dots, a_n = y$ 
    such that $a_{i-1} R a_i$ for every 
    $1 \leq i \leq n$.
    Then the transitive closure $\bar{R}$ of $R$ is
    anti-symmetric.
    \end{fact}
    \begin{proof}[Proof of Lemma \ref{lemma:Goulden_Yong_Order_On_T}]
    Every finite sequence of edges in $T$, with the property
    that every two consecutive edges $e$ and $f$ we have $e \prec_T f$, must form a path. 
    Now, between every two edges there is exactly one path, hence at most one sequence as above.
    Hence by Lemma \ref{fact:TransitiveClosureProperty}
    $<_T$ is anti-symmetric.
    It is clearly anti-reflexive, hence a strong order on the edges of $T$.
    \end{proof}
    
	We call $\le_T$ the {\em Goulden-Yong partial order} corresponding to $T$.
	
	\begin{observation}\label{conclusion:linear_extension}
	For every factorization $w = (t_1,\ldots,t_n) \in F_n$, the order
	$t_1 < t_2 < \ldots < t_n$ is a linear extension of the Goulden-Yong order $<_{G(w)}$.
	\end{observation}
	
	\begin{example}
	In Figure~\ref{fig:2}, the tree $T = G(w)$ yields the partial order satisfying $(1,4) <_T (4,6) <_T (4,5)$ and $(1,4) <_T (1,2) <_T (2,3)$.
	It is not a linear order. The order $(1,4) < (4,6) < (4,5) < (1, 2), (2,3)$
	is a linear extension of it.
	\end{example}

	\section{Convex caterpillars}
	
	In this section we prove Theorem \ref{thm:main} using the properties of convex caterpillars.
	\subsection{Basic properties of convex caterpillars}
	Let us use the following conventions.
	All  arithmetical operations on the elements of $[n]$ will be done modulo $n$. $[a, b]$ will denote the cyclic interval $\{a, a+1, \dots, b\}$. Using this notation, for edges $(a\ b), (a\ c)$ of a geometric non-crossing tree $T$, we have $(a\ c) <_T (a\ b)$ if and only if $b \in [a, c]$.
	
	The following result was proved in \cite{YK-CenterPaper}. We provide a somewhat different proof, the details of which will be used later.
	

	\begin{theorem}\cite[Theorem 3.2]{YK-CenterPaper}\label{thm:linear}
		The Goulden-Yong order on the edge set of a non-crossing geometric  tree $T$ is linear (total) if and only if $T$ is a convex caterpillar.
	\end{theorem}
	
	The following observation follows from the fact that a linear extension of a Goulden-Yong order $<_T$ on the edges of 
	geometric non-crossing tree $T$ corresponds to a factorization of the cycle 
	$(1\dots n)$ into $n-1$ transpositions.
	\begin{observation}\label{obs:consectuve_edges}
	If $T$ is a geometric non-crossing tree and $<_T$ is linear, then every two consecutive edges, viewed as transpositions in $\symm_n$, do not commute and therefore have a common vertex.
	\end{observation}
	
	The following lemma gives sufficient conditions for $<_T$ not to be linear.
	\begin{lemma}\label{lemma:non_linear}
	Let $T$ be a non-crossing geometric tree. In each of the following cases, the order $<_T$ is not linear.
	\begin{enumerate}
	    \item  There are edges $(a\ b), (c\ d), (e\ f)$ of $T$ such that $(a\ b) <_T (c\ d), (e\ f)$ and 
	    $c,d \in [a,b-1]$ and $e, f \in [b, a-1]$.
	    \item $T$ has edges $(a\ b), (c\ d), (e\ f)$ such that
	    $(c\ d), (e\ f) < (a\ b)$ and $c,d \in [b+1, a]$ and $e, f \in [a+1, b]$.
	\end{enumerate}
	\end{lemma}
	\begin{proof}
	We prove the first case, second one being similar by reversing directions. 
	Suppose that $<_T$ is linear and the first case holds. Note that for every $v \in [a+1, b-1]$ the edge $(v\ b)$ is smaller than $(a\ b)$ in $<_T$ because of counterclockwise relation of the edges, and the same is true for any edge $(a\ v)$ with $v \in [b+1\ a-1]$.
	Combining with non-crossing property of $T$ we find that any edge that is larger than $(a\ b)$ has either end-points in $[a, b-1]$ or $[b, a-1]$. Since it has both, there must be adjacent edges with endpoints in $[a, b-1]$ and $[b, a-1]$. However they are disjoint, hence commute, contradicting the fact that $<_T$ is linear order. 
	\end{proof}
	
	We are ready to prove Theorem \ref{thm:linear}.
	
	\begin{proof}[Proof of Theorem \ref{thm:linear}]
	If $T$ is a convex caterpillar. If the spine of $T$ is empty, then $T$ is a star, hence every two edges are comparable because they have a common vertex. Otherwise, let $(a\ a+1), (a+1\ a+2)\dots,(b-2\ b-1), (b-1\ b)$ be the spine of $T$. For every two edges $(k\ l)$ and $(k\ m)$ that share a common vertex $k$, $(k\ l) <_T (k\ m)$ if $m \in [k, l]$ where $[k\ l]$ denotes the cyclic interval $\{k,k+1,\dots,l \}$ where $l-k$ and $m-k$ have values between $1$ and $n-1$. Note that it is simply restatement of the fact that neighbors of $k$ are ordered counterclockwise.
	Hence, every two edges in the spine are comparable with $(a\ a+1) <_T (a+1\ a+2)\dots <_T (b-1\ b)$. It also implies that if $(k\ k+1)$ is an edge in the spine, then for every edge $(k\ l)$ that has $k$ as an end point, $(k\ l) <_T (k\ k+1)$ and for every edge $(k+1\ m)$ that has $k+1$ as an endpoint we have $(k\ k+1) <_T (k+1\ m)$. It follows that if $k$ and $m$ are endpoints of edges in the spine with $(k\ k+1),(k+1\ k+2)\dots,(m-1\ m)$, then for every edge $(k\ j)$ connected to $k$ and every edge $(m\ l)$ connected to $m$ we have $(k\ j) <_T (k\ k+1) <_T \dots (m-1\ m) \dots <_T (m\ l)$, hence every two edges that do not have common vertex are also comparable.
	
	To prove the converse statement, assume that $<_T$ has unique linear extension.
	Then $<_T$ is linear and we can sort the edges $(a_1\ b_1),\dots,(a_{n-1}, b_{n-1})$, and since every linear extension of $<_T$ corresponds to decomposition of the cycle $(1\ \dots \ n)$ into transpositions, we can view each edge as transposition.
	Next, note that since $<_T$ is linear, every two adjacent edges can not commute as transpositions, hence share a common vertex. 
	Now, note that the first edge must be of form $(i\ i+1)$ for some $i$. Assume that it $t_1 = (i\ j)$ where the length cyclic interval $[i j]$ is larger than $1$ and smaller than $n-1$. Since $T$ is a tree, there must exist a vertex $k$ in the cyclic interval $[i+1 j-1]$ and a vertex $m$ in the cyclic interval $[j+1 i-1]$ connected to either $i$ or $j$. Note that since every two consequent edges in $<_T$ must have a common vertex, $t_2$ must be connected to either $i$ or $j$. Assume without loss of generality that $t_2 = (j k)$ for some $k \in [j+1\ i-1]$. But every two consecutive edges in $<_T$ must have a common vertex, and every vertex adjacent to $t_2$ must have vertices in the interval $[j\ i-1]$ because of the non-crossing property of $<_T$. However, this implies that the first edge in in the interval $[i\ j-1]$ has no common vertex with the edge preceding it, which means that they commute as transpositions which contradicts the fact that $<_T$ is linear. 
	
	Now $i$ must be a leaf. For if we have an edge $(j\ i)$, $(i\ j) <_T (i\ i+1)$, contradicting the fact that $<_T$ is the first edge in $<_T$.  $m$ edges in $<_T$ the following hold:
	\begin{enumerate}
	    \item The end points of the first $m$ edges in $<_T$ form a cyclic interval $[j\ k]$.
	    \item The vertices $j, j+1, \dots, i-1, i$ are leaves in $T$.
	    \item The edges are $(i\ i+1), (i+1\ i+2), \dots (k-1\ k)$ are edges in $T$ and occur among the first $m$ edges.
	    \item Every edge that has $j, j+1 \dots, k-1$ as endpoint occurs among the first $m$ edges.
	    \item For the $m$-th edge in $<_T$, $t_m$ = $(k-1\ k)$ or $t_m = (k\ j)$. 
 	\end{enumerate}
  	Let $t_l$ denote the $l$-th edge in $<_T$.
 	The statement clearly holds for $m=1$. Assume that the statement holds for $m$. By induction hypothesis the $m$-th edge of $<_T$ is either $(k\ j)$ or $(k-1\ k)$ and linearity of $<_T$ and the induction hypothesis $t_{m+1}$ must have $k$ as an endpoint, because $j$ and $k$ can only be endpoints of the first $m$ edges by the hypothesis.
 	Next we show that $t_{m+1}$ is either $(k\ k+1)$ or $(k\ j-1)$. Assume $t_m = (k\ l)$ for $l\neq j, k+1$.  Then we have must have edges $(k\ l) <_T (s\ t), (u\ v)$ with $u,v \in [k\ l-1]$ and $s, t \in [l, k-1]$ contradicting Lemma \ref{lemma:non_linear}.
 	Now if $t_{m+1}=(k\ k+1)$, we are done, since the statements $1$ and $2$ hold by induction for $m$, $3$ and $5$ hold for $m+1$ and $4$ holds because $(k\ k+1)$ is the maximal edge in $<_T$ that has $k$ as an endpoint.
 	If $t_{m+1} = (k\ j-1)$, then for every $v \in [k+1, j-2]$ we have $(j-1\ v) <_T (j-1\ k)$ which is impossible, since $v\notin [j,k]$ contradicting the assumption.
 	On the other hand, for every $v \in [j\ k-1]$, $(j-1\ v)$ can not be a an edge, since by the assumption, since edges with endpoints $j, \dots k-1$ occur among the first $m$ edges. Hence, $j-1$ must be a leaf.
 	Again, it is easy to check that assumptions $1, 2, 3, 4, 5$ still hold for $m+1$.
 	Now if we substitute $m$ with $n-1$, we see that $T$ must be a geometric caterpillar, because by the construction, vertices that are not leaves are $i+1, i+2, \dots, k$ for some $k$, with edges $(i+1\ i+2), \dots, (k-1, k)$ connecting them.
	\end{proof}

	For example, the tree in Figure~\ref{fig:2} is a caterpillar, but not a convex one. The corresponding Goulden-Yong order is not linear.
	
	\begin{corollary}\label{cor:linear}
	A non-crossing geometric tree $T$ on $n$ vertices is a convex caterpillar if and only if there is a unique $w \in F_n$ such that $G(w) = T$. 
	\end{corollary}
	
	We shall henceforth identify a convex caterpillar $c \in \Ct_n$ with the corresponding sequence of transpositions $(t_1,\ldots,t_{n-1}) \in F_n$.

	\begin{prop}\label{prop:common_vertex}
		In a convex caterpillar $c=(t_1,\dots ,t_{n-1})$:
		\begin{enumerate}
		    \item 
		    Any two consecutive edges $t_i$ and $t_{i+1}$ share a common vertex.
		    \item 
		    The first edge $t_1$ is of the form $(a,a+1)$ for some $a$.
	        The same holds for the last edge $t_{n-1}$.
		\end{enumerate}
	\end{prop}
	\begin{proof}
	The first part of the proposition follows from the proof of \ref{thm:linear}.
    The second part is simply restatement of \ref{obs:consectuve_edges}.
	\end{proof}
	\begin{definition} Let $e$ be an edge of caterpillar $c$.
		\begin{enumerate}
			\item 
			We say that $e$ is a \emph{branch} if (at least) one of its endpoints is a leaf.
			\item 
			We say that $e$ is a \emph{link} if its endpoints have cyclically consecutive labels.
		\end{enumerate}
	\end{definition}
	
	By cautiously reading  the proof of theorem \ref{thm:linear}, we get the following observation.
	\begin{observation}\label{obs:first_last_edge}
		An edge of a convex caterpillar a $c$ is both a link and a branch if and only if it is either the first or the last edge of $c$.
	\end{observation}
	
 	\begin{lemma}\label{CaterpillarStructure}
		Let $c = (t_1,\dots, t_{n-1}) \in \Ct_n$. The following statements  hold. 
		\begin{enumerate}
			\item 
			The endpoints of the first $k$ edges form a cyclic interval in $[n]$, for every $1\leq k \leq n-1$.
			\item 
			If the first edge is $(i,i+1)$ then the endpoints of the first $k$ branches that are leaves are $i, i-1, \dots, i-k+1$, in that order.
			\item 
			If the first edge is $(i, i+1)$ then the first $k$ links are $(i,i+1),(i+1,i+2),\ \dots, \ (i+k-1, i+k)$.
			\item 
			The product of the first $k$ edges, viewed as transpositions, is equal to the cycle $(\ell, \ell+1, \dots, m)$ where $\ell$ is the leaf endpoint of the last branch among the first $k$ edges and $(m-1, m)$ is the last link among the first $k$ edges.
		\end{enumerate}
	\end{lemma}
	\begin{proof}
	Parts $1$, $2$ and $3$ follow from the proof of theorem \ref{thm:linear}. Part $4$ follows by induction and using the fact that that if the product of first $k$ edges (viewed as transpositions) is the cycle $(l\ l+1\ dots l+k)$ and where the cyclic interval is formed by the endpoints of first $k$ edges, then the $k+1$-th edge is either $(l+k\ l+k+1)$ or $(l+k\ l-1)$. Multiplying these we get the desired result.
	\end{proof}
	
	\begin{corollary}\label{conc:caterpillarstructure2}
	Every $c = (t_1,\dots, t_{n-1}) \in \Ct_n$ is completely determined by its first edge $t_1$ and the set of indices $i$ for which $t_i$ is a branch.
	\end{corollary}

		\subsection{A labeling of maximal chains}\label{sec:phi}
	The following labeling of maximal chains in the non-crossing partition lattice was introduced by Adin and Roichmain in \cite{AR_OMCNCPL} and is closely related to the the $EL$-labeling introduced by Bj\"orner in \cite{Bjorner}.
	In this section we describe this labeling, denoted by $\phi$. Its connection to the $EL$-Labeling of Bj\"orner will be discussed elsewhere.
	
	\medskip
	
	Recall, from Definition \ref{defn:Des-c}, the notion of descent set of a convex caterpillar. 
	
	Next, we show the connection to the descents defined by the map $\phi$ in \cite{AR_OMCNCPL}. First, let us describe $\phi$.
	For $w=(t_1,\dots,t_{n-1}) \in F_n$ define the partial products $\sigma_j = t_j\dots t_{n-1}$ with $\sigma_n = id$. By definition $\sigma_j = t_j\sigma_{j+1}$. For $1 \leq j \leq n-1$ define 
	$$A_j = \left\{1 \leq i \ n-1 : \sigma_j(i) > \sigma_{j+1}(i) \right\}.$$
	
	By the discussion preceding Definition 3.2 in \cite{AR_OMCNCPL}, we get the following statement.
	\begin{prop}\label{prop:phi_justification}
	The following hold.
	\begin{enumerate}
	    \item For each $1 \leq j \leq n-1$, $|A_j| = 1$.
	    \item The map $\pi_w$ defined by $$\pi_w(j) = i\ if\ A_j = \left\{ i \right\}$$ is a permutation in $\symm_{n-1}$.
	\end{enumerate}
	\end{prop}
	\begin{defn}\cite[Definition 3.2]{AR_OMCNCPL} Define $\phi:F_n\to\symm_{n-1}$ by $$\phi(w) = \pi_w.$$
	\end{defn}
	Define for each $w\in F_n:$ $$\Des(w)=\Des(\phi(w)).$$

	\subsection{Descents of convex caterpillars}
	
	We proceed to  calculate the restriction of $\phi$ to $Ct_n$.
	
	\begin{prop}\label{prop:alt_description_phi}
	Let $c\in Ct_n$ and let $\sigma_{j+1} = (\ k+1 \dots m)$. Then 
	$$\phi(c)(j) = 
	\begin{cases*}
	l & if $\sigma_{j+1}(n) \neq n$ \\
	l & if $\sigma_{j+1}(n) = n$ and $l < m$ \\
	m & if $\sigma_{j+1}(n) = n$ and $m < l$ \\
	\end{cases*}
	$$
	\end{prop}
	\begin{proof}
	By Lemma \ref{CaterpillarStructure} the product of the first $n-1-j$ transpositions is a cycle of form $(l\ l+1\ \dots m)$. However, this implies that  $\sigma_j$ is the cycle $(m\ m+1\ \dots\ l-1)$. Also, note that $(m-1, m)$ is the last link among the first $n-j$ edges of $c$.
	Hence $t_{j-1}$ equals either $(m-1\ m)$ or $(m\ l)$. 
	
	If $t_{j-1} = (m-1\ m)$ and $m-1 < m$, then 
	$$\sigma_{j-1}(m-1) = m > m-1 = \sigma_j(m-1)$$ 
	which implies that $\phi(c)(j-1) = m-1$. $m-1 > m$ then $m-1 = n$ and $$\sigma_{j-1}(l-1) = n > m = \sigma_j(l-1),$$
	hence, $\phi(c)(j-1) = l-1$.
	If $t_{j-1} = (m\ l)$ then if $l < m$, we have 
	$$\sigma_{j-1}(l) = m > l >\sigma_j(l)$$ and $\phi(c)(j-1) = l$ and if $m < l$ then $$\sigma_{j-1}(l-1) = l > m = \sigma_j(l-1)$$ and $\phi(c)(j-1) = l-1$. Note, that in all four cases, we get following combinatorial description of $\phi$ restricted to $Ct_n$.
	\end{proof}
	\begin{prop}\label{prop:phi_caterpillar_des}
	The descent set of a convex caterpillar, defined as in Definition~\ref{defn:Des-c}, 
	coincides with the descent set defined via the map $\phi$.
	\end{prop}
	\begin{proof}
	First, show that if $t_i = (a\ b)$ and $t_{i+1} = (b\ c)$ then $\phi(c)(i) > \phi(c)(i+1)$. Let $\sigma_{j+2} = (k\ k+1, \dots, m)$ such that $(k\ k+1)$ is the first link in $\sigma_{j+2}$ and $t_j = (k-1\ k)$ is the last link among the first $j-1$ edges. There are two possibilities. We have either $b = k$ or $c = k$ Suppose that $b = k$ holds then we have $t_j = (a\ k)$, $t_{j+1} = (c\ k)$ with $a > k$. Then by interval property of $\sigma_j$ of a caterpillar we have $b = m+1$ and $c = m+2$ with $m+2 > m+1$. By proposition \ref{prop:alt_description_phi}, $\phi(c)(j) = m+1,\ phi(c)(j+1) = m$ if $n\notin \left\{k,\dots,m\right\}$ and $\phi(c)(j) = m+2,\ \phi(c)(j+1)=m+1$. In both cases we have $\phi(c)(j) > \phi(c)(j+1)$, hence $j$ is a descent of $\phi(c)$. Second possibility is that $b \neq k$. In that case we have $(b\ c) = (k-1 k)$ and $(a\ b) = (m+1\ k-1)$ with $m+1 > k-1$. This implies that $n$ is not contained in the interval $(k-1\dots m+1)$ which means that $\phi(c)(j) = m$ and $\phi(c)(j+1) = k-1$ and again $j$ is a descent of $\phi(c)$.
	
	Now assume that $j$ is descent of $\phi$. 
	Let $\sigma_{j+2} = (k\ \dots \ m)$. Note that that there are four possibilites for $t_j, t_{j+1}$. 
	\begin{enumerate}
	    \item $t_j = (k-2\ k-1), t_{j+1} = (k+1\ k)$. In this case we have $\sigma_{j+1} = (k-1\ k\ \dots \ m)$, $\sigma_j = (k-2\ k-1\ \dots k)$.
	    By proposition \ref{prop:alt_description_phi}, either $\phi(c)(j+1) = k-1$ or $\phi(c)(j+1) = m$ if $k-1 = n$. We $\phi(j)(c) = m$ and or $k-2 = m$ if $k-2 = n$. Since we have $\phi(c)(j) > \phi(c)(j+1)$ we can not have $\phi(j) = k-2 > k-1 = \phi(j+1)$ because it would imply that $\phi(k-2) = n$ and this is not possible because $\phi(c)$ is permutation on $n-1$. Hence the possibilities that remain are either $\phi(j) = m > k-1 = \phi(j+1)$ or $\phi(j) = k-2 > \phi(j+1) = m$.
	    If $\phi(j) = m > k-1 = \phi(j+1)$ then we have $\sigma_j = (n\ k-1), \sigma_{j+1} = (k-1\ k)$ which means that $j$ is a descent of $c$. If 
	    $\phi(j) = k-2 > \phi(j+1) = m$, then $\sigma(j+1) = (n\ 1)$ and $\sigma(j) = (n-1\ n)$ which again implies that $j$ is a descent of $\phi$.
	    \item $t_j = (k-1\ k), t_{j+1} = (k\ m+1)$. Again by proposition \ref{prop:common_vertex} we have either $\phi(j) = k-1$ or $\phi(j) = m+1$ if $\sigma(j) = (n\ 1)$ and $\phi(j+1) = m$ or $\phi(j+1) = m+1$ if $\sigma_{j+2}(n) \neq n$. We must have either $\phi(c)(j) = m + 1 > m = \phi(c)(j+1)$. In this case we have $t_j = (n\ 1)$ and $t_{j+1} - (1\ m)$. Otherwise we have $\phi(j) = k-1 > m = \phi(j+1)$ or 
	    $\phi(j) = k-1 > m+1 = \phi(j+1)$. Both cases imply that $k-1 > m$ and thus $j$ is again descent of $c$.
	    \item $t_j = (k-1\ m+1), t_{j+1} = (k-1\ k)$. By proposition \ref{prop:common_vertex} we have either $\phi(c)(j) = m$ if $k-1 < m$ and $\phi(c)(j) = m+1$ if $k-1 > m, m+1$. We also have $\phi(c)(j+1) = k-1$ if $k-1 < k$ and $\phi(c)(j+1) = m$ if $t_{j+1} = (n\ 1) = (k-1\ k)$.
	    Clearly, the option $\phi(c)(j) = m+1 > k-1 > \phi(c)(j+1)$ is not 
	    possible because it implies $k+1 > m+1$ and $j$ is a descent of $\phi$.
	    Hence we have $\phi(c)(j) = m > k-1 = \phi(c)(j+1)$ which implies that $m+1 > k-1$. Hence we have $t_j = (k-1\ m+1), t_{j+1} = (k-1\ k)$ with $m+1 > k-1$ which implies that $j$ is a descent of $c$.
	    \item $t_j = (k\ m+2), t_{j+1} = (k\ m+1)$. If we have $m+1 > m+2$ then we have $m+1 = n$, hence $\phi(c)(j+1) = m = n$ and $\phi(c)(j) = 1$ which is not possible, since $j$ is a descent. Otherwise we have either $\phi(c)(j) = m+1, \phi(c)(j) = m$ or $\phi(c)(j) = m+2, \phi(c)(j+1) = m+1$ by proposition \ref{prop:alt_description_phi}. It is easy to check that in both cases $j$ is also a descent of $c$.
	    \end{enumerate}
	\end{proof}
	
	\subsection{Schur-positivity of convex caterpillars}
	
		\begin{definition}
		Let $c=(t_1,\dots ,t_{n-1})$ be a convex caterpillar and let $i$ be the index of the first edge that has $1$ as its endpoint. 
		The edge $t_i$ is called the \emph{main edge} of $c$ and the index $i$ is called the \emph{main index} of $c$, denoted $\mathtt{I}(c)$.
	\end{definition}
	
	For example, for
	$c = \left( (4,5),(5,6),(3,6),(1,6),(1,2) \right)$ we have $\mathtt{I}(c) = 4$.
	
	\medskip
	
	Using  Lemma~\ref{CaterpillarStructure} we prove the following explicit description of the descents of a convex caterpillar $c$, based on $\mathtt{I}(c)$ and on the geometry of $c$.
	
	\begin{lemma}\label{DescentBranchLemma1}
		Let $c \in \Ct_n$ and $i \in [n-2]$. Then:
		\begin{enumerate}
		  	\item 
		  	For $1 \le i < \mathtt{I}(c)-1$, $i$ is a descent of $c$ if and only if $t_{i+1}$ is a branch of $c$. 
		    \item 
		    For $i = \mathtt{I}(c)-1$,  $i$ is always a descent of $c$.
			\item 
			For $i=\mathtt{I}(c)$, $i$ is a descent of $c$ if and only if $1$ is not a leaf of $c$.
			\item 
			For $\mathtt{I}(c) < i \le n-2$, $i$ is a descent of $c$ if and only if $t_i$ is a branch of $c$. 
		\end{enumerate}
	\end{lemma}
	\begin{proof} We prove each case separately
		\begin{enumerate} 
			\item First suppose that $i = \mathtt{I}(c)$. Then $t_i = (a\ b)$ and $t_{i+1} = (b\ 1)$ for some $2 \leq a, b \leq n$.
			Obviously, $a > 1$ and therefore $t_i$ is a descent. 
			
			\item If $t_{i+1}$ is a branch, then $t_i = (a\ b),\ t_{i+1} = (a\ c)$ for some $a, b, c > 1$. If $(a\ b)$. By lemma \ref{CaterpillarStructure} $b = a-1$ if $(a\ b)$ is a link or $b = c+1$ if $(a\ b)$ is a branch, the endpoints of the first $i+1$ edges form the cyclic interval $[c, a]$. Since $i < i(c)$ we $1 < c < b < a$, therefore $i$ is a descent. On the other hand if $t_{i+1}$ is a link then $t_i = (a, b)$, $t_{i+1} = (b\ b+1)$ and $a$ is between $b+1$ and $b$ in $<_b$. Because $1 < a$, we have $a < b$, thus $i$ is not a descent.
			\item Now if $i = \mathtt{I}(c)$ and $1$ is a leaf. Then we have 
			$ t_i = (a\ 1),\ t_{i+1} = (a\ b).$
			Obviously $1 < b$ and $t_i$ is not a descent. In contrast, if $t_i$ is a link, then $t_i = (1\ a),\ t_{i+1}=(1\ b)$ and since $a$ and $b$ are sorted counterclockwise and are both greater than $1$ in the cyclic order $<_1$ we have $b < a$, thus $i$ is a descent.  
			\item \label{partC} Now suppose that $i > \mathtt{I}(c)$. Then if $t_i$ is a branch we have 
			$t_i = (a\ b)\ t_{i+1} = (a\ c)$ where $b$ and $c$ are ordered counterclockwise and $a < b, c < n$ which implies that $c < b$ and that $t_i$ is a descent. On the other hand, if $t_i$ is a link we have  $t_i = (a\ a+1),\ t_{i+1} = (a+1\ k)$ where $k > a$, and $i$ is not a descent.
		\end{enumerate}
	\end{proof}

	Combining  Lemmas~\ref{CaterpillarStructure}  and~\ref{DescentBranchLemma1} and Corollary~\ref{conc:caterpillarstructure2},
	we deduce the following key proposition.
	
	\begin{prop}\label{UniqueLemma}
		A convex caterpillar $c$ is determined uniquely by $\mathtt{I}(c)$ and $\Des(c)$.
	\end{prop}
	\begin{proof}
		By Lemma \ref{conc:caterpillarstructure2}, it suffices to show that the pair $(\mathtt{I}(c),\Des(c))$  
		determines 
		the first edge and branches. Note that by observation \ref{obs:first_last_edge}, first and last edges are always branches.
		
		Denote $k:= \mathtt{I}(c)$.
		For $i < k$, combine Lemmas~\ref{DescentBranchLemma1} and~\ref{CaterpillarStructure} to determine whether the $i$-th edge is a branch or link. By Part 2 of Lemma~\ref{DescentBranchLemma1} we know whether $1$ is a leaf or not and whether $e_k$ is a branch or not. In both cases, applying Parts $2$ and $3$ of Lemma~\ref{CaterpillarStructure}, we determine the first edge. The branches with indices larger than $k$ are determined by Part $3$ of Lemma~\ref{DescentBranchLemma1}. Hence the first edge and the branches are completely determined by the descent set and the $k=\mathtt{I}(c)$ as desired.
	\end{proof}
	
	The next lemma describes the possible values of $\mathtt{I}(c)$, given the descent set of $c$. 
	
	\begin{lemma}
		Let $c \in \Ct_n$. Then either $\mathtt{I}(c) = 1$ or $\mathtt{I}(c) - 1 \in \Des(c)$.
	\end{lemma}
	\begin{proof}
		Let $X\subseteq\left[n-2\right]$ and suppose that $\Des(c) = X$. We show that $\mathtt{I}(c) = 1$ or $\mathtt{I}(c)\in X + 1$. 
		It is clear that if $\mathtt{I}(c)\neq 1$ then there exists a $i \in \mathtt{des}_C(1)+1$ such that $\mathtt{I}(c)=i$ by Part $1$ of Lemma \ref{DescentBranchLemma1}.
	\end{proof}
	
	\begin{lemma}
		For every subset $J\subseteq\left[n-2\right]$ and every $i\in\left(1+J\right)\cup\left\{1\right\}$, there exists a unique $c \in \Ct_n$ such that $\Des(c)=J$ and $\mathtt{I}(c)=i$. 
	\end{lemma}
	\begin{proof}
		Recall that every caterpillar is determined by its first edge and the true branches, where every $J\subseteq\left\{2,\dots,n-2\right\}$ can appear as the set of the true branches of a caterpillar.
		Placing $\mathtt{I}(c)$ after $x\in J$ results in proper set of true branches, which in turns defines a caterpillar. Now, suppose that $i\notin J$. Then $(1, 2)$ can be first edge of the leaf, since $1$ is a leaf, hence $1$ is not a descent, and branches correspond to the members of $X$. If $1\in X$, then $(n\ 1)$ can be first edge, with the rest of branches defined by the descents.
	\end{proof}
	
	\begin{corollary}\label{conclusion}
		For every subset $J\subseteq\left[n-2\right]$, the number of convex caterpillars with descent set $J$ is equal to $\left| J \right| + 1$.
	\end{corollary}
	
	The following observation is well known.

	
	\begin{observation}\label{lemma:hook}
	For every $0 \le k \le n-1$
	\[
	\{\Des(T) \,:\, T\in \SYT(n-k,1^k)\} = \{J \subseteq [n-1] \,:\, |J|=k\},
	\]
	each set being obtained exactly once.
	\end{observation}
	
	\begin{proof}[Proof of Theorem~\ref{thm:main}]
	Combine Corollary~\ref{conclusion} with Observation~\ref{lemma:hook} and Theorem~\ref{G1} to deduce
	\[
	\Q(Ct_n)=\sum_{k=0}^{n+1}(k+1)\sum\limits_{\substack{J\subseteq [n-1] \\ |J|=k}} \F_{n,J} =
	\sum_{k=0}^{n+1}(k+1)s_{n-k,1^k}.
	\]
	\end{proof}
	
	\smallskip
	
	\noindent{\bf Acnowledgements.} 
	This work forms part of a PhD research conducted under the supervision of Professors Ron Adin and Yuval Roichman.

%
	

\end{document}